\documentclass{article}

\newcommand{\Bbb}{\bf}

\newcommand{\kk}{{\bf k}}
\newcommand{\nn}{{\bf n}}
\newcommand{\rr}{{\bf r}}
\newcommand{\ZZ}{{\Bbb Z}}

\newcommand{\RR}{{\cal R}}
\newcommand{\OO}{{\bf O}}
\newcommand{\II}{{\cal I}}
\newcommand{\JJ}{{\cal J}}
\newcommand{\tJ}{\widetilde{\JJ}}

\newcommand{\lam}{\lambda}
\newcommand{\sig}{\sigma}
\newcommand{\Om}{\Omega}
\newcommand{\taum}{ \tau^{ \mbox{\tiny max} } }
\newcommand{\taur}{\tau^{\rr}}

\newcommand{\eqto}{\stackrel{\rm \sim}{\to}}
\newcommand{\zi}{(\zeta^*)^{-1}}

\newcommand{\Fl}{\mbox{\rm Fl}}

\newcommand{\Stab}{\mbox{\rm Stab}}
\newcommand{\rank}{\mbox{\rm rank}}
\newcommand{\codim}{\mbox{\rm codim}}
\newcommand{\mod}{\mbox{\rm mod}}
\newcommand{\Block}{\mbox{\rm Block}}

\begin{document}

\centerline{\large \bf DEGENERACY SCHEMES}
\centerline{\large \bf AND SCHUBERT VARIETIES}
\vspace{1em}
\centerline{\sc \large v.~lakshmibai and peter magyar
\footnote
{Both authors partially supported by the
National Science Foundation.}
}
\vspace{1em}
\centerline{29 August 1997}
\vspace{2em}

\noindent
{\small {\sc abstract.}  A result of Zelevinsky states that
an orbit closure in the space of representations of the
equioriented quiver of type $A_h$ is in bijection with the
opposite cell in a Schubert variety of a partial flag
variety $SL(n)/Q$.  We prove that Zelevinsky's bijection is
a scheme-theoretic isomorphism,
which shows that the universal
degeneracy schemes of Fulton are reduced and
Cohen-Macaulay in arbitrary characteristic.
}
\\[2em]
Among all algebraic spaces, the best understood are the
flag varieties and their Schubert subvarieties.
They first appear as interesting examples, but acquire a
general importance in the theory of characteristic
classes of vector bundles.

Fulton has recently given a theory \cite{F} of
>``universal degeneracy loci'', characteristic
classes associated to maps among vector bundles,
in which the role of Schubert varieties is taken by
certain degeneracy schemes.  The underlying varieties
of these schemes
arise in the theory of quivers: they are the
closures of orbits in the
space of representations of the equioriented quiver
$A_h$.  The variety of complexes is another example
of this class.

By a remarkable but little-known result of Zelevinsky
\cite{Z} (c.f.~\cite{L}), the
above quiver varieties can be identified
set-theoretically with open subsets
of Schubert varieties of a partial flag variety.
In this paper, we prove a scheme-theoretic strengthening
of Zelevinksy's identification:
the ``naive'' determinantal conditions defining each
quiver variety as a set generate the same ideal as
the Plucker equations defining
(the opposite cell of) the corresponding Schubert variety.
Since the latter ideal is well understood via
Standard Monomial Theory, we conclude that
the quiver schemes defined by the determinantal
equations are reduced and their singularities are identical
to those of Schubert varieties.  In particular,
the quiver varieties are
Cohen-Macaulay, answering a question of Fulton.

Our results extend those of Abeasis, Del Fra, and Kraft
\cite{A}, \cite{AF}, \cite{AFK},
who showed that the quiver varieties are Cohen-Macaulay
with rational singularities over a field of characteristic zero,
and that the determinantal conditions generate
the reduced ideals of the quiver varieties of codimension one.
Also, our methods are similar to those of Gonciulea and Lakshmibai
\cite{GL}.
\\[1em]
{\bf Contents.}
{\bf 1.} Zelevinsky's bijection \quad
 1.1 Quiver  varieties \quad
 1.2  Schubert varieties \quad
 1.3  The bijection $\zeta$ \quad
 1.4  The actions of $B$, $Q$, and $G_{\nn}$ \quad
{\bf 2.} Plucker coordinates and determinantal ideals \quad
 2.1 Coordinates on the opposite big cell \quad
 2.2 The main theorem \quad
 2.3 Proof of the main theorem: determinant identities \quad
 2.4 Degeneracy schemes

\section{Zelevinsky's bijection}
\subsection{Quiver  varieties}

Fix an $h$-tuple of non-negative integers
$\nn = (n_1,\ldots,n_h)$
and a list of vector spaces $V_1,\ldots, V_h$
over an arbitrary field $\kk$
with respective dimensions $n_1,\ldots,n_h$.
Define the {\it variety of quiver representations}
(of dimension $\nn$, of the equioriented quiver of type
$A_h$) to be the affine space $Z$ of all
$(h\!-\!1)$-tuples of linear maps $(f_1,\ldots,f_{h\!-\!1}):$
$$
V_1 \stackrel{f_1}{\to} V_2
 \stackrel{f_2}{\to} \cdots \stackrel{f_{h\!-\!2}}{\to} V_{h\!-\!1}
 \stackrel{f_{h\!-\!1}}{\to} V_h \ .
$$
If we endow each $V_i$ with a basis, we get $V_i \cong \kk^{n_i}$
and
$$
Z \cong M(n_2 \times n_1) \times \cdots
\times M(n_{h} \times n_{h\!-\!1}) ,
$$
where $M(k\times l)$ denotes the affine space of matrices
over $\kk$ with $k$ rows and $l$ columns.
The group
$$
G_{\nn} = GL(n_1) \times \cdots \times GL(n_h)
$$
acts on $Z$ by
$$
(g_1,g_2,\cdots,g_h) \cdot (f_1,f_2,\cdots,f_{h\!-\!1})
= (g_2 f_1 g_1^{-1}, g_3 f_2 g_2^{-1},\cdots,
g_{h}f_{h\!-\!1} g_{h\!-\!1}^{-1}),
$$
corresponding to change of basis in the $V_i$.

Now, let $\rr = (r_{ij})_{1 \leq i \leq j \leq h}$
be an array of non-negative integers with $r_{ii} = n_i$,
and define $r_{ij} = 0$ for any indices other than
$1\leq i\leq j \leq h$.  Define
$$
Z^{\circ}(\rr) = \{(f_1,\cdots,f_{h\!-\!1}) \in Z
\ \mid\  \forall\, i\!<\!j,\ \rank (f_{j\!-\!1} \cdots f_i : V_i \to V_j)
 = r_{ij} \}.
$$
(This set might be empty for a bad choice of $\rr$.)
\\[1em]
{\bf Proposition.} {\it The $G_{\nn}$-orbits of
$Z$ are exactly the sets $Z^{\circ}(\rr)$
for $\rr=(r_{ij})$ with
$$
r_{i,j\!-\!1} -
r_{i,j} -
r_{i\!-\!1,j\!-\!1} +
r_{i\!-\!1,j}  \geq 0,\quad
\forall\ 1\! \leq\! i\! <\! j\! \leq\! h.
$$
}

\noindent {\bf Proof.} This is a standard result of algebraic
quiver theory \cite{BGP}, \cite{G}, \cite{W}.
Since this theory is not well known among geometers,
we recall it here.

Consider the abelian category
$\RR$ of quiver representations
whose objects are sequences of linear maps
$(V_1 \stackrel{f_1}{\to} \cdots \stackrel{f_{h\!-\!1}}{\to} V_h)$,
where the $V_i$ are {\it any} vector spaces of
{\it arbitrary} dimension.
A morphism of $\RR$ from the object
$(V_1 \stackrel{f_1}{\to} \cdots \stackrel{f_{h\!-\!1}}{\to} V_h)$
to the object
$(V'_1 \stackrel{f'_1}{\to} \cdots \stackrel{f'_{h\!-\!1}}{\to} V'_h)$
is defined to be an $h$-tuple of linear maps $(\phi_i:V_i \to V'_i)$
such that each square
$$
\begin{array}{ccc}
V_i & \stackrel{f_i}{\to} & V_{i+1} \\
\mbox{\tiny $\phi_i$} \downarrow &&
 \downarrow \mbox{\tiny $\phi_{i+1}$} \\
V'_i & \stackrel{f'_i}{\to} & V'_{i+1}
\end{array}
$$
commutes.

Direct sum of objects is defined componentwise, and it is
known (Krull-Remak-Schmidt Theorem) that any object $R \in
\RR$ can be written uniquely as a direct sum of
the indecomposable objects
$$
\begin{array}{ccccc}
R_{ij} =
(0 \to \cdots \to 0\to &
\!\!\!\!\kk\!\!\!\! & \eqto \cdots \eqto \kk \to &
\!\!\!\!0\!\!\!\! & \to \cdots \to 0) \\
&\!\!\!\! V_i \!\!\!\!&&
\!\! \!\!V_j\!\! \!\! &
\end{array}
$$
for $1 \leq i<j \leq h+1$ (corresponding to the
positive roots of the root system $A_h$).
That is,
$$
R \cong \bigoplus_{1 \leq i<j \leq h+1} m_{ij} R_{ij}
$$
for unique multiplicities $m_{ij} \in \ZZ^+$.

Our variety $Z$ consists of representations
with fixed $(V_i)$ and all possible $(f_i)$.
Two points of $Z$ are in the same $G_{\nn}$-orbit
exactly if they are isomorphic as objects in $\RR$.
So the orbits correspond to arrays
 $(m_{ij})_{1 \leq i<j \leq h+1}$ with $m_{ij} \in \ZZ^+$
and $n_i = \sum_{k\leq i< l} m_{kl}$.

We can compute the rank numbers $(r_{ij})$ from
the multiplicities $(m_{ij})$:
$$
r_{ij} = \sum_{k\leq i<j<l} m_{kl},
$$
and conversely
$$
m_{ij}=
r_{i,j\!-\!1} -
r_{i,j} -
r_{i\!-\!1,j\!-\!1} +
r_{i\!-\!1,j} .
$$
Hence the arrays $(r_{ij})$ with the stated conditions
classify the $G_{\nn}$-orbits on $Z$. $\bullet$
\\[1em]
We define the {\it quiver variety}
$$
Z(\rr)
=\{(f_1,\cdots,f_{h\!-\!1}) \in Z
\mid  \forall i,j,\ \rank (f_{j\!-\!1} \cdots f_i : V_i \to V_j)
\leq r_{ij}\}.
$$
Finally, we have the dimension formula due
to Abeasis and Del Fra \cite{AF}.
\\[1em]
{\bf Propsition.}
$$
\dim Z(\rr) = \dim G_{\nn} -\!\!\!
\sum_{1\leq i\leq j \leq h}
(r_{ij}-r_{i,j+1}) (r_{ij}-r_{i-1,j}).
$$

\subsection{Schubert varieties}

Given $\nn=(n_1,\cdots,n_h)$, for $1 \leq i \leq h$ let
$$
a_i = n_1 + n_2 + \cdots +n_i
\qquad \mbox{and} \qquad
n = n_1  + \cdots + n_h \ .
$$
For positive integers $i \leq j$, we shall frequently use
the notations
$$
[i,j] = \{ i, i+1, \ldots, j\}, \qquad\qquad [i] = [1,i]\ .
$$

Let $\kk^n$ be a vector space (over our arbitrary field $\kk$)
with standard basis $e_1,\ldots,e_n$.  Consider its general
linear group $GL(n)$, the subgroup $B$ of upper-triangular
matrices, and the parabolic subgroup $Q$ of block upper-triangular
matrices
$$
Q = \{ (a_{ij} \in GL(n) \mid a_{ij}=0 \
\mbox{whenever}\ j\leq a_k <i
\ \mbox{for some}\ k \}\ .
$$
A {\it partial flag of type $(a_1<a_2<\cdots <a_h=n)$ }
(or simply a {\it flag}) is a sequence of supspaces
$U. = (U_1 \subset U_2 \subset \cdots \subset U_h = \kk^n)$
with $\dim U_i = a_i$.
Let $E_i = \langle e_1,\ldots,e_{a_i}\rangle$
the span of the first $a_i$ coordinate vectors, and
$E'_i = \langle e_{a_i+1},\ldots,e_n\rangle$ the natural
complementary subspace to $E_i$, so that
$E_i \oplus E'_i = \kk^n$.
Call $E. = (E_1 \subset E_2 \subset \cdots)$
the {\it standard flag}.
Let $\Fl$ denote the set of all flags $U.$ as above.

$\Fl$ has a transitive $GL(n)$-action induced from
$\kk^n$, and $Q = \Stab_{GL(n)}( E.)$, so we may identify
$\Fl \cong GL(n)/Q$, \ $g\!\! \cdot\!\! E. \leftrightarrow gQ$\ .
The {\it Schubert varieties} are the closures of $B$-orbits
on $\Fl$.  Such orbits are usually indexed by certain
permutations of $[n]$, but we prefer to use
{\it flags of subsets} of $[n]$, of the form
$$
\tau = (\tau_1 \subset \tau_2 \subset\cdots \subset \tau_h = [n]),
\qquad \#\tau_i=a_i\ .
$$
(A permutation $w: [n]\to[n]$ corresponds
corresponds to the subset-flag with
$\tau_i = w[a_i] = \{w(1),w(2),\ldots,w(a_i)\}$.
This gives a one-to-one correspondence between cosets
of the symmetric group $S_n$ modulo the Young subgroup
$S_{n_1} \times \cdots \times S_{n_h}$, and subset-flags.)

Given such $\tau$, let
$E_i(\tau) = \langle e_j \mid j \in \tau_i \rangle$
be a coordinate subspace of $\kk^n$, and
$E.(\tau) = (E_1(\tau) \subset E_2(\tau) \subset \cdots) \in \Fl$.
Then we may define the {\it Schubert cell}
$$
\begin{array}{rcl}
X^{\circ}(\tau) &= &B\cdot E(\tau)\\
&=& \left\{(U_1\subset U_2\subset\cdots)\in \Fl\ \ \left|\
\begin{array}{c}
\dim U_i \cap \kk^j = \#\, \tau_i \cap [j]\\[.2em]
1\leq i \leq h,\ 1\leq j \leq n
\end{array}
\right.\right\}
\end{array}
$$
and the {\it Schubert variety}
$$
\begin{array}{rcl}
X(\tau) &= &\overline{X^{\circ}(\tau)}\\
&=& \left\{(U_1\subset U_2\subset\cdots)\in \Fl\ \ \left|\
\begin{array}{c}
\dim U_i \cap \kk^j \geq \#\, \tau_i \cap [j]\\[.2em]
1\leq i \leq h,\ 1\leq j \leq n
\end{array}
\right.\right\}
\end{array}
$$
where $\kk^j = \langle e_1,\ldots,e_j\rangle \subset \kk^n$.

We define the {\it opposite cell} $\OO \subset \Fl$
to be the set of flags in general position with respect
to the spaces
$E'_1 \supset \cdots \supset E'_{h-1}$:
$$
\OO = \{(U_1\subset U_2\subset\cdots)\in \Fl\ \mid\
U_i \cap E'_{i}=0\}.
$$
We also define $Y(\tau) = X(\tau) \cap \OO$, an open subset of
>$X(\tau)$.  By abuse of language, we call $Y(\tau)$ the
{\it opposite cell} of $X(\tau)$, even though it is not a cell.

\subsection{The bijection $\zeta$}

We define a special subset-flag
$\taum = (\taum_1 \subset \cdots \subset\taum_h = [n])$
corresponding
to $\nn = (n_1,\ldots,n_h)$.
We want $\taum_i$ to contain numbers as large as possible
given the constraint $[a_{i\!-\!1}]\subset \taum_i$.
Namely, we define $\taum_i$ recursively by
$$
\taum_h = [n];\quad \taum_{i} = [a_{i\!-\!1}]
\cup \{ \mbox{largest $n_i$ elements of $\taum_{i+1}$}\}.
$$
Furthermore, given $\rr = (r_{ij})_{1\leq i\leq j\leq h}$
indexing a quiver variety, define a subset-flag $\taur$ to
contain numbers as large as possible given the
constraints
$$
\#\, [a_j]\cap \taur_i =
\left\{ \begin{array}{cl}
a_i -r_{i,j+1} & \mbox{for}\ i\leq j \\
a_j& \mbox{for}\ i> j \\
\end{array} \right.
$$
Namely,
$$
\taur_i = \{\,
\underbrace{1\ldots a_{i\!-\!1}}_
{\mbox{\small $a_{i\!-\!1}$}}
\ \underbrace{. \ldots\ldots a_{i}}_
{\mbox{\small $r_{ii}\!-\!r_{i,i+1}$}}
\ \underbrace{.\ldots\ldots a_{i+1}}_
{\mbox{\small $r_{i,i+1}\!-\!r_{i,i+2}$}}
\ \underbrace{.\ldots\ldots a_{i+2}}_
{\mbox{\small $r_{i,i+2}\!-\!r_{i,i+3}$}}\ \ldots\
\ \underbrace{.\ldots\ldots n_{\mbox{}}}_
{\mbox{\small $r_{i,h}$}}
\}
$$
where we use the visual notation
$$
\underbrace{\cdots\cdots a}_{\mbox{\small $b$}} = [a-b+1,a].
$$
Note that $r_{ij} -r_{i,j+1} \leq n_j$,
so that each $\taur_i$ is a list of increasing integers,
and that $r_{ij}-r_{i,j+1}\leq r_{i+1,j}-r_{i+1,j+1}$,
so that $\taur_i \subset \taur_{i+1}$.  Thus $\taum$ and
$\taur$ are indeed subset-flags.

Now define the Zelevinsky map
$$
\begin{array}{rccc}
\zeta: & Z & \to &\Fl \\
& (f_1,\ldots,f_{h\!-\!1}) & \mapsto &(U_1\subset U_2 \subset \cdots)
\end{array}
$$
where
$$
U_i = \{ (u_1,\ldots,u_h)\in
\kk^{n_1}\!\oplus\! \cdots \!\oplus\! \kk^{n_h} = \kk^n \mid
\forall\, j>i,\ u_{j+1} = f_j(u_j)\}.
$$
In terms of coordinates, if we identify the linear maps
$(f_1,\ldots,f_{h\!-\!1})$ with the matrices $(A_1,\ldots,A_{h\!-\!1})$,
and identify $\Fl \cong GL(n)/Q$,
we have
$$
\zeta(A_1,\ldots,A_{h-1}) =
\left(
\begin{array}{ccccc}
I_1 & 0 & 0 & 0& \cdots \\
A_1 & I_2 &0 & 0& \cdots \\
A_2 A_1 & A_2 & I_3 & 0& \cdots \\
A_3 A_2 A_1 & A_3 A_2 &A_3 & I_4 & \cdots  \\[-.4em]
\vdots & \vdots & \vdots &  \vdots&
\end{array}
\right) \ \
\mod \ \ Q
$$
where $I_i$ is an identity matrix of size $n_i$.
\\[1em]
{\bf Theorem.}  {\it (Zelevinsky \cite{Z})\\
(i) $\zeta$ is a bijection of $Z$ onto its image $Y(\taum)$:
\quad $\zeta : Z \eqto Y(\taum)$.\\
Also, \\[-1em]
$$
\mbox{}\hspace{-2em} (*)\qquad
Y(\taum)=
\{(U_1\subset U_2 \subset \cdots) \ \mid \
\forall\ i,\ \ E_{i-1} \subset U_i,
\ \ U_i \cap E'_{i} = 0 \}.
$$
(ii) $\zeta$ restricts to a bijection from $Z(\rr)$
onto $Y(\taur)$:\quad $\zeta : Z(\rr) \eqto Y(\taur)$.\\
Also, \\[-1em]
$$
\mbox{} \hspace{-2em} (**) \quad
Y(\taur)=\left\{(U_1\subset U_2 \subset \cdots) \ \left| \
\begin{array}{c}
 \forall\ i\leq j,\quad \dim \, E_j \cap U_i \geq a_i-r_{i,j\!+\!1}
,\\[.3em]
 E_{i-1} \subset U_i,
\quad U_i \cap E'_{i} = 0\end{array}
\right.\right\}\ .
$$
}

\noindent{\bf Proof.}  Obviously $\zeta$ is injective.
To prove (i), we first show that $\zeta(Z)$ is equal to the
right hand side of equation $(*)$.  One inclusion is clear.

To show the other inclusion,
consider any $U.$ with $E_{i\!-\!1} \subset U_i$
and $U_i\cap E'_{i}=0$ for all $i$.
Let $\pi_i:\kk^n = E_i \oplus E'_i \to E_i$ be the projection.
Then $\pi_{h\!-\!1}$
restricts to an isomorphism $U_{h\!-\!1} \eqto E_{h\!-\!1}$,
so there exists an inverse linear map
$$
\mbox{id} \oplus f_{h\!-\!1} :
E_{h\!-\!1} \to E_{h\!-\!1}\oplus\kk^{n_h}
$$
such that
$$
U_{h\!-\!1} = \mbox{Graph}(f_{h\!-\!1})
 \subset E_{h\!-\!1} \oplus \kk^{n_h}
=\kk^n\ .
$$
Since $E_{h\!-\!2}\subset U_{h\!-\!1}$, we have
 $f_{h\!-\!1}(E_{h\!-\!2})=0$.
Next, $\pi_{h\!-\!2}$ restricts to an isomorphism
$U_{h-2} \eqto E_{h-2}$, and
there exists a linear map
$\tilde{f}_{h\!-\!2}: E_{h\!-\!2} \to E'_{h\!-\!2}$ with
$\tilde{f}_{h\!-\!2}(E_{h\!-\!3})=0$
such that
$$
U_{h\!-\!2} = \mbox{Graph}(\tilde{f}_{h\!-\!2})
 \subset E_{h\!-\!2} \oplus E'_{h\!-\!2}
=\kk^{n}.
$$
Since $U_{h\!-\!2}\subset U_{h-1}$, we have
$$
\tilde{f}_{h\!-\!2} = (f_{h\!-\!2},\, f_{h\!-\!1}f_{h\!-\!2})
$$
for some $f_{h\!-\!2}:E_{h\!-\!2} \to \kk^{n_{h\!-\!1}}$.
Continuing in this way, we find that $U.\in \zeta(Z)$.


Thus it suffices
to show that $(*)$ is valid.  Again, the inclusion
$\subset$ is clear.
Now consider a  flag $U.$ satisfying
$E_{i\!-\!1} \subset U_i$ for all $i$.
Then we will show that $U.$ must
satisfy $\dim(\kk^i \cap U_j) \geq \#\, [i] \cap \taum_j$
for all $1\leq i\leq n$,\, $1 \leq j\leq	h$.
Acting by $B$ does not change $\dim U_i \cap \kk^j$,
so we may assume our $U.$
is a flag of coordinate subspaces $U. = E(\mu)$ for some
$\mu=(\mu_1\subset\cdots\subset \mu_{h}=[n])$
with $[a_{i\!-\!1}]\subset \mu_i$ for all $i$,
so that $\dim U_i \cap \kk^j = \# \mu_i \cap [j]$.
Then by the definition of $\taum$, we must have
$$\forall\, j, \qquad
\# \mu_{h\!-\!1} \cap [j] \geq
\# \taum_{h\!-\!1}\cap [j],\qquad
\# \mu_{h\!-\!2} \cap [j] \geq
\# \taum_{h\!-\!2}\cap [j],\
\ldots\ldots
$$
This proves $(*)$, and hence part (i).

The proof of (ii) is similar.
Clearly $Y(\taur) \subset Y(\taum) = \zeta(Z)$.
For any flag $U.=\zeta(f_1,\ldots,f_{h\!-\!1})$, we have
$$
\begin{array}{rcl}
\dim\, E_j \cap U_i &=&
\dim E_{i\!-\!1} + \dim \mbox{Ker}(f_j f_{j\!-\!1} \cdots f_i)\\
&=& \dim E_{i\!-\!1} + \dim V_i - \rank( f_j f_{\!j-\!1}\cdots f_i)\\
&=& a_i -\rank( f_j f_{j\!-\!1}\cdots f_i).
\end{array}
$$
Hence $\dim\, E_{j} \cap U_i
\geq a_i - r_{i,j\!+\!1}$ if and only if $U. \in \zeta(Z(\rr))$,
so that $\zeta(Z(\rr))$ is equal to the right hand side
of $(**)$.   But the conditions on the right side of $(**)$
are enough to force the flag $U.$ to lie in the Schubert
variety $X(\taur)$ on the left hand side, as in part (i).
$\bullet$

\subsection{The actions of $B$, $Q$ and $G_{\nn}$}

Let $W = S_n$ and
$W_{\nn} = S_{n_1} \times \cdots \times S_{n_h}$ a Young
subgroup.  Let $W_{\nn}$ act on the
the coset space $W/W_{\nn}$ by left multiplication.
Then we may consider $\taum$ as a coset
in $W/W_{\nn}$ which is Bruhat-maximal within its
$W_{\nn}$ orbit.
Since $W_{\nn}$ is the Weyl group of $Q$,
this means that the $B$-action on
the Schubert variety $X(\taum)$
extends to a $Q$-action.

We may embed $G_{\nn}$ into $Q$ as the block diagonal
matrices, so that $G_{\nn}$ acts on $X(\taum)$ and
in fact on the open subvariety $Y(\taum)$.
Then $\zeta: Z \to Y(\taum)$ is equivariant
with respect to the $G_{\nn}$-action.
\\[1em]
Now we relate our combinatorial formalism to
that in Zelevinsky's original paper \cite{Z}.
We have just seen that our $\taur$ correspond
to certain double cosets in $W_{\nn}\backslash W / W_{\nn}$.
Following Zelevinsky, we may index such double cosets
by {\it block permutation matrices}, which are defined
to be the $h\times h$ arrays $T=(t_{ij})$ of non-negative integers
with row and column sums equal to the $n_i$ , so that
for all $1\leq i,j\leq h$,
$$
\sum_{i=1}^h t_{ij} = n_j \qquad \sum_{j=1}^h t_{ij} = n_i\ .
$$
(If all $n_i=1$, this defines an ordinary permutation matrix.)

A permutation $w \in W$ correponds to the block permutation
matrix $\Block(w)$ defined by partitioning the ordinary
permutation matrix of $w$ into blocks, and summing all entries
in each block:
$$
\Block(w) = (t_{ij}) \qquad
t_{ij} = \#\ [a_{i\!-\!1}+1, a_i] \cap w[a_{j\!-\!1}+1,a_j]\ .
$$
The block map induces a one-to-one correspondence between
double cosets $W_{\nn}\backslash W / W_{\nn}$ and block permutation
matrices.

Zelevinsky's map takes $Z(\rr)$ to $Y(\taur)$ for each
$\rr=(r_{ij})$.  Recall from the proof of Proposition 1.1
that the rank numbers $r_{ij}$, $1\leq i\leq j \leq h$,
can be computed from certain multiplicities $m_{ij}$,
$1\leq i<j \leq h+1$.
Then the block permutation matrix corresponding to $\taur$
is given by
$$
\left(
\begin{array}{ccccc}
m_{12} & m^*_{1} & 0 & 0 &\cdots \\
m_{13} & m_{23} & m^*_2 & 0 & \cdots \\
m_{14} & m_{24} & m_{34} & m^*_3 & \cdots \\
\vdots & \vdots & \vdots & \vdots &
\end{array}
\right)
$$
where
$$
m^*_{i} = \sum_{k < i+1 < l} m_{kl}\ .
$$

\section{Plucker coordinates and determinantal ideals}

For a variety $X$ embedded in an affine space $V$
over an infinite field $\kk$, the {\it vanishing ideal} $\II$
of $X$ is the set of polynomial functions on $V$ which restrict
to zero on $X$.  However, if $\kk$ is a finite field,
we modify this definition in the usual way:
the vanishing ideal is the the set of polynomials on
$V$ which are zero on the points of $X$ over the algebraic
closure of $\kk$:
$$
\II=\{f \in \kk[V] \mid f(x)=0\ \, \forall\,
x\in X(\overline{\kk})\}.
$$
The ideal $\II$ is necessarily reduced (radical).

\subsection{Coordinates on the opposite big cell}

Consider the opposite
cell $\OO \subset GL(n)/Q$.  It is easily seen that
$\OO$ consists of those cosets which have a unique
representative $A$ of the form
$$
A = (a_{kl}) =
\left( \begin{array}{ccccc}
I_1& 0 & 0& \cdots & 0\\
A_{21} & I_2 & 0 & \cdots &0\\
A_{31}& A_{32} & I_3 &\cdots & 0\\
\vdots&\vdots&\vdots &&\vdots\\
A_{h1}& A_{h2}& A_{h3}&\cdots& I_{h}
\end{array} \right)\  \mod \   Q,
$$
where $I_i$ is the identity matrix of size $n_i$,
and $A_{ij}$ is an arbitrary matrix of size $n_i \times n_j$.
That is, $\OO$ is an affine space with coordinates
$a_{kl}$ for those positions $(k,l)$
with $1 \leq l \leq a_i <k \leq n$ for some $i$.
Its coordinate ring is the polynomial ring
$$
\kk[\OO] = \kk[a_{kl}].
$$

For a matrix $M \in M(k\times l)$ and subsets
$\lam\subset [k]$, $\mu \subset [l]$, let
$\det M_{\lam\times \mu}$ be the minor with row indices $\lam$
and column indices $\mu$.
Now let $\sig \subset [n]$ be a subset of size
$\#\sig = a_i$ for some $i$.  Define the {\it Plucker
coordinate} $p_{\sig} \in \kk[\OO]$ to be the
$a_i$-minor of our matrix $A$ with row indices $\sig$ and
column indices the interval $[a_i]$:
$$
 p_{\sig}=p_{\sig}(A)=\det A_{\sig \times [a_i]}.
$$
Define a partial order on Plucker coordinates by:
$$
\sig \leq \sig'
\quad \Longleftrightarrow \quad
\begin{array}{c}
\sig = \{\sig(1)<\sig(2)<\cdots<\sig(a_i)\},\\
\sig' = \{\sig'(1)<\sig'(2)<\cdots<\sig'(a_i)\},\\
\sig(1)\leq \sig'(1),\ \sig(2) \leq \sig'(2),
\cdots, \sig(a_i) \leq \sig'(a_i).
\end{array}
$$
This is a version of the Bruhat order.
\\[1em]
{\bf Proposition.} {\it
Let $\tau = (\tau_1 \subset \cdots \subset \tau_{h} = [n])$
be a subset-flag and $Y(\tau)$ the intersection of the
Schubert variety $X(\tau)$ with the opposite cell $\OO$.  Then
the vanishing ideal $\II(\tau) \subset \kk[\OO]$ of
$Y (\tau) \subset \OO$ is generated by those Plucker coordinates
$p_{\sig}$ which are incomparable with one of the $p_{\tau_i}$:
$$
\II(\tau) = \langle p_{\sig} \mid \exists\, i,\ \#\sig = a_i,
\ \sig \not\leq \tau_i \rangle.
$$
}

\noindent
{\bf Proof.}  This follows from well-known results of
Lakshmibai-Musili-Seshadri in
Standard Monomial Theory (see e.g.~\cite{MS},\cite{LS}).

\subsection{The main theorem}

Denote a generic element of the quiver space
$ Z = M(n_2\times n_1) \times \cdots
\times M(n_{h}\times n_{h\!-\!1})$
by $(A_1,\ldots,A_{h-1})$, so that the coordinate ring
of $Z$ is the polynomial ring in the entries of all the matrices
$A_i$.  Let $\rr = (r_{ij})$ index the quiver variety
$Z(\rr) = \{(A_1,\ldots,A_{h-1}) \mid
\rank\, A_{j-1}\cdots A_i \leq r_{ij}\}$.

Let $\JJ(\rr) \subset \kk[Z]$ be the ideal generated by
the determinantal conditions implied by the definition
of $Z(\rr)$:
$$
\JJ(\rr) = \left\langle \det(A_{j-1} A_{j-2} \cdots A_i)_
{\lam\times\mu}
\ \left| \
\begin{array}{c}
i\leq j,\ \lam \subset [n_j],\ \mu \subset [n_i] \\[.2em]
\#\lam = \#\mu = r_{ij}+1
\end{array}
\right.
\right\rangle\ .
$$
Clearly $\JJ(\rr)$ defines $Z(\rr)$ set-theoretically.
\\[1em]
{\bf Theorem.} {\it
$\JJ(\rr)$ is a prime ideal and is the vanishing ideal
of $Z(\rr)\subset Z$.  There are isomorphisms of
reduced schemes
$$
Z(\rr) = \mbox{Spec}(\kk[Z]\,/\,\JJ(\rr)) \cong
\mbox{Spec}(\kk[\OO]\,/\,\II(\taur))
= Y(\taur).
$$
That is, the quiver scheme $Z(\rr)$ defined by $\JJ(\rr)$ is
equal to the reduced, irreducible variety $Y(\taur)$,
the opposite cell of a Schubert variety.
}
\vspace{.5em}

\noindent {\bf Proof.} Consider the map of \S1.3,
$\zeta: Z \eqto Y(\taum) \subset \OO$.  It is clear that
$\zeta$ is an algebraic isomorphism onto its image, since
it is injective on points and on tangent vectors.
(In fact, in appropriate coordinates
$\OO \cong Z\times V$ for some affine space $V$, and
for a certain polynomial function $\phi:Z\to V$, \, $\zeta$
is equivalent to the map
$Z \to Z\times V$, \ $z \mapsto (z,\phi(z))$.\,)
Thus by Proposition 2.1 we have the exact sequence
$$
0\to \II(\taum) \to \kk[\OO]
\stackrel{\zeta^*}{\to}
\kk[Z] \to 0\ .
$$

Let $\tJ(\rr) \subset \kk[Z]$ be the (reduced)
vanishing ideal of $Z(\rr) \subset Z$.
Clearly $\JJ \subset \tJ$. Since
$\zeta$ maps $Z(\rr)$ isomorphically onto
$Y(\taur)$ by Theorem 1.3,  we have
$\zi\tJ(\rr) = \II(\taur)$ by Proposition 2.1. Hence
$$
\begin{array}{rcl}
Z(\rr) = \mbox{Spec}(\kk[Z]\,/\,\tJ(\rr)) &\cong&
\mbox{Spec}(\kk[\OO]\,/\,\zi\tJ(\taur))\\
&=&\mbox{Spec}(\kk[\OO]\,/\,\II(\taur))
\ = \ Y(\taur).
\end{array}
$$
Furthermore, $\JJ(\rr) = \tJ(\rr)$ if and only if
$\zi\JJ(\rr) = \zi\tJ(\rr)$; and $\JJ(\rr)$ is prime
if and only if $\zi\JJ(\rr)$ is prime.
Thus to show the Theorem, it suffices to prove
$$
\zi\JJ(\rr) = \II(\taur).
$$
But clearly $\zi \JJ(\rr) \subset \zi\tJ(\rr) = \II(\taur)$,
so we are left with the opposite inclusion
$$
\zi \JJ(\rr) \supset \II(\taur),
$$
which we will prove in the next section.

\subsection{Proof of the main theorem: determinant identities}

We define ideals $\II_0, \II_1, \II_2 \subset \kk[\OO]$
generated by certain minors of the generic matrix $A \in \OO$
at the end of \S1.2:
$$
\II_0 =\zi\JJ(\rr) \hspace{3.5in} \mbox{}
$$
\vspace{-1.8em}
$$
= \II(\taum) +
\left\langle
\det (A_{j,j\!-\!1} A_{j\!-\!1,j\!-\!2}\!\cdots\!
A_{i\!+\!1,i} )_{\lam\times\mu}
\left|
\begin{array}{c}
i\!<\!j,\ \
\lam\! \subset\! [n_j],\  \ \mu\!\subset\! [n_i]\\[.2em]
 \# \lam = \#\mu = r_{ij}+1
\end{array}
\right. \right\rangle
$$
\vspace{.4em}
$$
\II_1\ =\ \II(\taum) +
\left\langle \det A_{\lam\times\mu}\ \left|\
\begin{array}{c}
i<j,\ \
\lam \subset [a_j\!+\!1,n],\ \ \mu\subset [a_i]\\[.2em]
\# \lam = \#\mu = r_{ij}\!+\!1
\end{array}
\right. \right\rangle
$$
\vspace{.4em}
$$
\II_2\ =\  \II(\taur)\ =\
\left\langle \det A_{\sig\times [a_i]}\ \left|\
\begin{array}{c}
1\leq i \leq h\!-\!1,\ \ \sig \subset [n]\\[.2em]
\# \sig = a_i,\ \ \sig \not \leq \taur_i
\end{array}
\right. \right\rangle
$$
To finish the proof of Theorem 2.2, we will show
$$
\II_0 \supset \II_1 \supset \II_2\ .
$$

\noindent
{\bf Lemma 1.} {\it
Let $X = (x_{ij})$ and $Y= (y_{kl})$ be matrices of
variables $x_{ij}$, $y_{kl}$ generating a polynomial ring.
Let $\JJ_{X}$  (resp.~$\JJ_{Y}$) be the ideal generated by all
$r\!+\!1$-minors of $X$ (resp.~$Y$).  Then $\JJ_X$ and $\JJ_Y$
both contain all $r\!+\!1$-minors of the product $XY$.
}
\vspace{.5em}

\noindent {\bf Proof.}
$$
\det (XY)_{\lam\times\mu}
= \sum_{\nu} \det X_{\lam\times \nu}\, \det Y_{\nu\times\mu}.
\quad \bullet
$$

\noindent {\bf Lemma 2.} {\it
Let $(A_1,\ldots,A_{h-1})$ be a generic element of $Z$, and
for $i \leq j$
let $\JJ_{ij}$ be the ideal generated by
all $r\!+\!1$-minors of the $n_j \times n_i$ product matrix
$A_j\cdots A_i$.
Then $\JJ_{ij}$ contains all $r\!+\!1$-minors of the
$(n\!-\!a_{j\!-\!1}) \times a_i$ matrix
$$
\widetilde{A} =
\left( \begin{array}{cccc}
A_j\!\! \cdots\! A_{1} & A_j\!\!\cdots\! A_2& \cdots& A_j\!\! \cdots\! A_i \\
A_{j\!+\!1}\!\! \cdots\! A_{1} & A_{j\!+\!1}\!\!\cdots \!A_2& \cdots&
A_{j\!+\!1}\!\!
\cdots\! A_i \\
\vdots&\vdots&&\vdots\\
A_h\!\! \cdots\! A_{1} & A_h\!\!\cdots\! A_2& \cdots& A_h\!\! \cdots A_i\! \\
\end{array} \right)
$$
}
\vspace{.5em}

\noindent {\bf Proof.} Note that we can factor the matrix
$$
\widetilde{A} = \left(\!\!\!\begin{array}{c}
I_{j} \\ A_{j\!+\!1} \\ \vdots \\ A_h\! \cdots\! A_{j\!+\!1}
\end{array} \!\!\!\right)
 \cdot\, A_j\!\!\cdots\!\! A_i\, \cdot\
( A_{i\!-\!1}\!\!\cdots\!\! A_1,\
A_{i\!-\!1}\!\!\cdots\!\! A_2,\ \cdots\ ,\ A_{i\!-\!1},\ I_i)
$$
Now apply Lemma 1 twice.
\\[1em]
{\bf Lemma 3.}\qquad
$ \II_0 \supset \II_1\ .$
\\[.5em]
{\bf Proof.}  Let $\lam \subset [a_j\!+\!1]$,\
$\mu \subset [a_i]$, $\#\lam = \#\mu = r_{ij}+1$.
Then clearly
$$
\det A_{\lam\times \mu}
\in \zi( \det \widetilde{A}_{\lam\times\mu} )\ .
$$
%
%
Hence by Lemma 2, the generators of $\II_1$ lie in $\II_0$.
$\bullet$
\\[1em]
{\bf Lemma 4.} {\it (Gonciulea-Lakshmibai)\
Let $A$ be a generic element of $\OO$.
Let $1 \leq t\leq a_i$, \ $1\leq s \leq n$, and
$\sig = \{ \sig(1)<\sig(2)<\cdots<\sig(a_i)\} \subset [n]$
with $\sig(a_i-t+1) \geq s$.  Then $p_{\sig}(A)$
belongs to the ideal of $\kk[\OO]$ generated by $t$-minors
of $A$ with row indices $\geq s$ and column indices $\leq a_i$.
}
\vspace{.5em}

\noindent{\bf Proof.}  Choose $\sig'\subset [s,n] \cap \sig$
with $\#\sig' =t$, and let $\sig'' = \sig \setminus \sig'$.
Then the Laplace expansion of $p_{\sig}(A)$ with respect
to the rows $\sig'$, $\sig''$, gives
$$
p_{\sig}(A) = \det A_{\sig\times[a_i]}
= \sum_{\lam' \cup \lam''= [a_i]}
\!\!\!\pm \det A_{\sig'\times \lam'}
\det A_{\sig''\times\lam''},
$$
where the sum is over all partitions of the interval $[a_i]$.
The first factor of each term in the sum
is of the form required. $\bullet$
\\[1em]
{\bf Lemma 5.}\qquad  $\II_1 \supset \II_2\ .$
\\[.5em]
{\bf Proof.}
Let $\sig \subset [n]$ with $\# \sig = a_i$,\
$\sig \not \leq \taur_i$ for some $i$,\ $1 \leq i\leq h\!-\!1$.
Now, $\taur_i$ has the largest
possible entries such that
$$
\taur_i(a_i-r_{i,j+1})\leq a_j,\qquad \forall \,j\geq i ,
$$
so $\sig \not \leq \taur_i$ must violate this condition for
some $j$:
$$
\sig(a_i-r_{i,j+1}) \geq a_j+1,\qquad \exists\, j\geq i.
$$
Hence by Lemma 4, $p_{\sig}(A)$ is in $\II_1$. $\bullet$
\\[1em]
The Main Theorem  2.2 is therefore proved.

\subsection{Degeneracy schemes}

Fulton \cite{F} defines
the universal degeneracy scheme
$\Om_w$ associated to a permutation $w \in S_{m+1}$
as follows.
Fix $2m$ vector spaces $F_1,F_2,\ldots,F_m,E_m,\ldots,E_2,E_1$
with $\dim F_i = \dim E_i = 1$,
and let
$$
Z = M_{2\times 1}\times M_{3\times 2} \times
\cdots \times M_{m\times m\!-\!1} \times M_{m\times m}
\times M_{m\!-\!1\times m} \times\cdots\times M_{1\times 2}
$$
be the quiver space of all maps of the form
$$
F_1 \to F_2 \to \cdots \to F_m \to E_m \to \cdots \to E_2 \to E_1.
$$
(For convenience we will refer to these maps and their compositions
by symbols such as $F_i\to F_j$ and $F_i\to E_j$.)
Define rank numbers
$$
r(F_i,E_j) = \# \,[i] \cap w[j], \quad 1\leq i,j\leq m
$$
$$
r(F_i,F_j) = i\quad r(E_j,E_i)=i\quad 1\leq i<j\leq m
$$
and let $\rr_w$ be the array of these numbers.
Then let
$$
\Om_w = Z(\rr_w) \subset Z,
$$
the variety of all quiver representations satisfying
$$
\rank(F_i\to E_j) \leq \#\, [i]\cap w[j], \quad 1\leq i,j\leq m
$$
$$
\rank(F_i\to F_j) \leq i,\quad
\rank(E_j\to E_i)\leq i,\quad 1\leq i<j\leq m.
$$
(The latter conditions are clearly superfluous.)
More precisely, define $\Om_w$ as a scheme
by the same determinantal equations defining $Z(\rr_w)$
in \S2.2.
\\[1em]
{\bf Proposition.} {\it The scheme $\Om_w$ over an arbitrary
field $\kk$ is reduced and is isomorphic to
the opposite cell of a Schubert variety in
$\Fl = GL(n)/Q$, a partial flag variety of
>$\kk^n$, where $n = 2(1+\cdots+m) = m(m+1)$.
In particular,
$\Om_w$ is irreducible, Cohen-Macaulay, and normal,
and has rational singularities.}
\\[1em]
{\bf Proof.} This follows since Schubert varieties
are known to have these properties (see e.g.~\cite{R}).
$\bullet$
\\[1em]
{\bf Proposition.}
$$
\codim_Z\, \Om_{w} = \ell(w),
$$ where
$$\ell(w)=\#\left\{ (i,j) \ \left|\
\begin{array}{c} 1\leq i,j\leq m\!+\!1 \\[.4em]
i<j,\quad w(i) > w(j)
\end{array}\right.\right\}
$$
is the Bruhat length.
\\[1em]
{\bf Proof.}  By the dimension formula of
Abeasis and Del Fra \cite{AF} (c.f.~\S1.1 above), we have
$$
\begin{array}{rcl}
\dim \Om_w & = & \dim G_{\nn} -
\sum_{1\leq i,j\leq m}
(r(F_i,E_j)-r(F_i,E_{j\!-\!1}))
(r(F_i,E_j)-r(F_{i\!-\!1},E_{j})) \\[.5em]
&&\quad
-\sum_{1 \leq i\leq j\leq m} (r(F_i,F_j)-r(F_i,F_{j\!+\!1}))
(r(F_i,F_j)-r(F_{i\!-\!1},F_{j})) \\[.5em]
&&\quad
-\sum_{1 \leq i\leq j\leq m} (r(E_j,E_i)-r(E_j,E_{i\!-\!1}))
(r(E_j,E_i)-r(E_{j\!+\!1},E_{i}))
\\[.5em]
&&\quad -(r(F_m,F_m)-r(F_m,E_m))(r(F_m,F_m)-r(F_{m-1},F_m))
\\[.9em]
&=& 2(1^2 + 2^2 +\cdots+m^2)\\[.5em]
&&-
\sum_{1\leq i,j\leq m} \#\,([i]\cap w[j])\setminus ([i]\cap w[j\!\!-\!\!1])
\cdot
\#\,([i]\cap w[j])\setminus ([i\!\!-\!\!1]\cap w[j]) \\[.9em]
&=& 2(1^2 + 2^2 +\cdots+m^2) \ -\
\#\left\{\ (i,j) \ \left|\
\begin{array}{c} 1\leq i,j\leq m +1 \\[.4em]
w^{-1}(i) \leq j,\quad i \geq w(j)\
\end{array}\right.\right\}
\\[1em]
&=&  2(1^2 + 2^2 +\cdots+m^2) - m - \ell(w)\ .
\end{array}
$$
On the other hand,
$$
\dim Z = 2(1\cdot 2 + 2 \cdot 3 + \cdots + (m\!-\!1) \cdot m) + m^2.
$$
Hence
$$
\begin{array}{rcl}
\codim_Z\, \Om_w &=&
2(1\cdot 2 + \cdots + (m\!-\!1) \cdot m) + m^2
- 2(1^2 + \cdots m^2) + m + \ell(w)\\[.3em]
&= & \ell(w).
\end{array}
$$
\vspace{-2em}

\mbox{}$\hfill \bullet$
\\[1em]
{\bf Concluding remarks}
\\[.5em]
Denote $Z=Z(m)$ and $\Om_w= \Om_w(m)$ to emphasize the
dependence on $m$.  Consider $S_{m+1} \subset S_{m+2}$
in the usual way.  Then there is a natural map
$\pi: Z(m\!+\!1) \to Z(m)$ given by forgetting the middle
two spaces, and we may easily see the stability property:
$$
\Om_w(m\!+\!1) = \pi^{-1} \Om_w(m).
$$
The map $\pi$ is a fiber bundle over some open set of $\Om_w(m)$, and
>from the previous proposition, the generic
fibers of $\pi: \Om_w(m\!+\!1)\to \Om_w(m)$ have the same
dimension as the generic fibers of $\pi: Z(m+1) \to Z(m)$.

Finally, we note that $\Om_w(m\!+\!1)$ is closely related to
a Schubert variety of $\Fl' = GL(m\!+\!1)/B$,
the complete flag variety of $\kk^{m\!+\!1}$,
a much smaller flag variety
than that of Zelevinsky's bijection (c.f.~Fulton \cite{F} \S 3).
Namely, consider the open set $Z^{\circ}(m\!+\!1)$
of elements $F_1 \to \cdots\to  E_1$
with $F_{i}\to F_{i\!+\!1}$ injective,
$E_{i\!+\!1}\to E_i$ surjective, and $F_{m\!+\!1}\to E_{m\!+\!1}$
bijective.
Then we have a principal fiber bundle
$$
\begin{array}{cccc}
\psi: & Z^{\circ} & \to & \Fl' \times \Fl' \\
& (F_1 \to \cdots \to E_1) & \mapsto & (V.\, ,U.)
\end{array}
$$
where
$$
V_i = \mbox{Im}(F_i \to E_{m+1}) \qquad
U_i = \mbox{Ker}(E_{m+1} \to E_{m+1-i})\ .
$$
Now, letting $\Om_w^{\circ}(m\!+\!1) =
\Om_w(m\!+\!1) \cap Z^{\circ}(m\!+\!1)$,
an open subet of $\Om_w(m\!+\!1)$, we have
$$
\Om_w^{\circ}(m\!+\!1) = \psi^{-1}\
\{\,  (V.,U.) \in \Fl' \times \Fl'\  \mid \
V_i \cap U_j \leq \#\, w w_0 [i] \cap [j]\, \}.
$$
where $w_0$ is the longest element of $S_{m\!+\!1}$.
It is well-known that the subset of $\Fl' \times \Fl'$
on the right is a fiber bundle
over $\Fl'$ with fiber equal to the Schubert variety
$X(w w_0) \subset \Fl'$.

\end{document}